\documentclass{amsart}
\usepackage{amssymb} 
\usepackage{graphicx}      
\usepackage[all]{xy}
\newcommand{\CC}{{\mathbb C}} 
\newcommand{\ZZ}{{\mathbb Z}} 
\newcommand{\NN}{{\mathbb N}}

\newcommand{\liea}[1]{\mathfrak{#1}} 
\newcommand{\lieg}[1]{\mathrm{#1}}

\newcommand{\rk}{\operatorname{rk}} 
\newcommand{\ad}{\operatorname{ad}}

\newcommand{\supp}{\operatorname{supp}}
 
\theoremstyle{plain} 
\newtheorem{thm}{Theorem}[section] 
 
\newtheorem{lm}[thm]{Lemma} 
\newtheorem{cor}[thm]{Corollary} 
\newtheorem{conj}{Conjecture} 
 
\theoremstyle{definition} 
\newtheorem{re}[thm]{Remark} 
\newtheorem{ex}{Example} 
\newtheorem{defn}{Definition} 
 
\begin{document} 
 
\title{Richardson elements for classical Lie algebras} 
\author{Karin Baur} 
\thanks{Supported by Freie Akademische Stiftung and 
by a DARPA Grant}
\address{Karin Baur, Department of Mathematics, University of 
California, San Diego, USA}
\email{kbaur@math.ucsd.edu}
 
\date{January 20, 2005} 
 
\begin{abstract}
Parabolic subalgebras of semi-simple Lie algebras 
decompose as $\liea{p}=\liea{m}\oplus\liea{n}$ 
where $\liea{m}$ is a Levi factor and 
$\liea{n}$ the corresponding nilradical. 
By Richardsons 
theorem \cite{ri}, 
there exists an open orbit under the action of 
the adjoint group $P$ on the nilradical. 
The elements of this dense orbits are known as 
Richardson elements. 

In this paper we describe a normal form 
for Richardson elements in the classical case. 
This 
generalizes a construction for 
$\liea{gl}_N$ of Br\"ustle, Hille, Ringel and 
R\"ohrle \cite{bhrr} to the other classical 
Lie algebra and it extends the authors normal forms 
of Richardson 
elements for nice parabolic subalgebras of simple 
Lie algebras to arbitrary parabolic subalgebras 
of the classical Lie algebras \cite{b04}.
As applications we obtain a description of 
the support of Richardson elements and we recover 
the Bala-Carter label of the orbit of Richardson 
elements. 
\end{abstract}

\maketitle

%
\section*{Introduction} 
%
%
The goal of this paper is to describe Richardson 
elements for parabolic subalgebras of the 
classical Lie algebras. 

Let $\liea{p}$ be a parabolic subalgebra of 
a semi-simple Lie algebra $\liea{g}$ over $\CC$ 
and 
$\liea{p}=\liea{m}\oplus\liea{n}$ a Levi 
decomposition. By a fundamental theorem 
of Richardson \cite{ri} there always exist 
elements $x$ in the nilradical $\liea{n}$ 
such that $[\liea{p},x]=\liea{n}$. In other 
words, if $P$ is the adjoint groups 
of $\liea{p}$, 
then the orbit $P\cdot x$ is dense in $\liea{n}$. 
It is usually called the Richardson orbit. 
Richardson orbits have been studied for 
a long time and there are many open questions 
related to this setting. 
Our goal is to give explicit representatives 
for Richardson elements. In the case of 
$\liea{gl}_n$ there is a beautiful way to 
construct Richardson elements that has 
been described by Br\"ustle, Hille, Ringel 
and R\"ohrle in~\cite{bhrr}. Furthermore, 
Richardson elements with support in the 
first graded part $\liea{g}_1$ (where the 
grading is induced from the parabolic 
subalgebra) have been given for all simple 
Lie algebras in~\cite{b04}. 

However, these constructions do not work in 
general for classical Lie algebras. 
To fill this gap, we have modified the 
existing approaches to obtain Richardson elements 
for parabolic subalgebras of the 
classical Lie algebras. We do this using 
certain simple line diagrams. They correspond 
to nilpotent matrices with at most one 
non-zero entry in each row and in each 
column. 
We show that for most parabolic subalgebras, there exists 
a simple line diagram that defines 
a Richardson element. 
But there are cases where this is not possible 
as we will see. 
We expect that the representatives we describe 
will give more insight and hopefully answer 
some of the open questions. 
One of the interesting questions in the theory of 
Richardson elements is the structure of 
the support of a Richardson element. 
Recall that any parabolic subalgebra $\liea{p}$ 
induces a $\ZZ$-grading of $\liea{g}$, 
\[
\liea{g}=\oplus_{i\in\ZZ}\liea{g}_i\quad\text{with} 
\quad\liea{p}=\oplus_{i\ge 0}\liea{g}_i= 
\liea{g}_0\oplus(\bigoplus_{i>0}\liea{g}_i)
\]
where $\liea{g}_0$ is a Levi factor and 
$\liea{n}:=\oplus_{i>0}\liea{g}_i$ the corresponding 
nilradical. For details, we refer to our joint 
work with Wallach,~\cite{bw}. 
The support of a Richardson element 
$X=\sum_{\alpha\text{ root of }\liea{n}}k_{\alpha}X_{\alpha}$ 
are the roots of the nilradical $\liea{n}$ 
with $k_{\alpha}\neq 0$ (where $X_{\alpha}$ spans the 
root subspace $\liea{g}_{\alpha}$). 
The support $\supp(X)$ of $X$ lies in the subspace 
$\liea{g}_1\oplus\dots\oplus\liea{g}_k$ 
for some $k\ge 1$. 
For the normal form of Richardson elements 
we can determine the minimal $k_0$ such 
that $\supp(X)\subset$ 
$\liea{g}_1\oplus\dots\oplus\liea{g}_{k_0}$. 
We also recover the Bala-Carter label  of the 
dense orbit of Richardson elements, also 
called the {\itshape type} of the orbit. 
The Bala-Carter label is used in the classification 
of nilpotent orbits of simple Lie algebras, given 
in~\cite{bc}. 
For a description of these labels 
see chapter 8 of~\cite{cm}. 
The type of any nilpotent orbit in a classical 
Lie algebra has been described by 
Panyushev \cite{pan} 
in terms of the partitions of the orbit.  

\noindent
Before we describe our results and explain 
the structure of this 
article, we need to fix some notation. 
If $\liea{p}$ is a parabolic subalgebra 
of a semi-simple Lie algebra $\liea{g}$ 
we can assume that $\liea{p}$ 
contains a fixed Borel subalgebra. In this 
case we say that $\liea{p}$ is standard. 
If $\liea{m}$ is a Levi factor of $\liea{p}$ 
we say that $\liea{m}$ is standard if it 
contains a fixed Cartan subalgebra $\liea{h}$ 
that is contained in the fixed Borel 
subalgebra.

From now on we will assume that $\liea{g}$ is 
a classical Lie algebra, unless stated otherwise. 
As usual, the Cartan subalgebra 
consists of the diagonal matrices and the 
fixed Borel subalgebra is the set of upper 
triangular matrices. 
Then a standard Levi factor has the shape of 
a sequence of square matrices (blocks) on 
the diagonal and zeroes outside. 
In the case of $\liea{so}_{2n}$, we have 
to be careful: we will only consider 
parabolic subalgebras where $\alpha_n$ and 
$\alpha_{n-1}$ are both roots of the Levi 
factor or both roots of the nilradical 
or $\alpha_{n-1}$ a root of the Levi 
factor and $\alpha_n$ a root of the nilradical. 
In other words the case $\alpha_n$ a root of the Levi 
factor and $\alpha_{n-1}$ a root of the 
nilradical will be identified with this last 
case since the two parabolic subalgebras 
are isomorphic. 
So our standard 
$\liea{p}$ or $\liea{m}$ are uniquely 
defined by 
the sequence $d:=\underline{d}=(d_1,\dots,d_r)$ 
of the sizes of these blocks (and by specifying 
the type of the Lie algebra). 

We start by defining line diagrams for 
dimension vectors in section~\ref{se:line-diag}. 
It will turn out that each horizontal line diagram 
corresponds uniquely to 
elements of the nilradical of the 
parabolic subalgebra of $\liea{sl}_n$ of the given 
dimension vector. 
In section~\ref{se:rich-theory}
we gather the necessary properties 
of Richardson elements. 
In section~\ref{se:sl-case} we show that 
horizontal line diagrams in fact correspond 
to Richardson elements of the given parabolic 
subalgebra. The construction of such diagrams 
for $\liea{gl}_n$ appears first in~\cite{bhrr}. 
We have alreday mentioned that 
for the other classical Lie algebras, the 
horizontal line diagrams do not give 
Richardson elements. 
In general, the matrix obtained 
is not an element of the Lie 
algebra in question. Thus we will introduce 
generalized line diagrams in section~\ref{se:BCD-type} 
to obtain Richardson elements for parabolic 
subalgebras of the symplectic and orthogonal 
Lie algebras. As a by-product we obtain 
the partition of a Richardson element for the 
so-called simple parabolic subalgebras. 
The last section discusses the cases 
where line diagrams do not produce 
Richardson elements. For these we will 
allow ``branched'' diagrams. 
In the appendix we add examples illustrating 
branched diagrams. 
%
%
\section{Line diagrams}\label{se:line-diag}
%
%
Let $d=(d_1,\dots,d_r)$ be a dimension 
vector, i.e. a sequence 
of positive integers. 
Arrange $r$ columns of $d_i$ dots, top-adjusted. 
A {\it (filled) line diagram} for $d$, 
denoted by $L(d)$, 
is a collection of lines joining vertices of 
different columns such that each vertex is 
connected to at most one vertex of a column left 
of it and to at most one vertex of a column right 
of it 
and such that it cannot be extended by any line. 

We say that it is a {\it (filled) horizontal line diagram} 
if all edges are horizontal lines. Such a diagram 
will be denoted by $L_h(d)$. 
We will always assume that the 
line diagrams are filled and omit the term `filled'. 
Line diagrams are not unique. However, for each dimension 
vector there is a unique horizontal line diagram.

\begin{ex}
As an example, consider the dimension vector $(3,1,2,3)$ 
and 
three line diagrams for it, the last one horizontal. 
$$
{\small
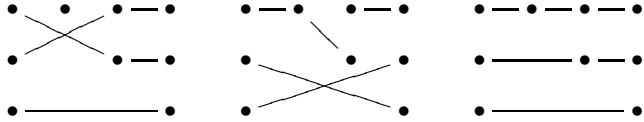
\xymatrix@-5mm{
\bullet\ar@{-}[rrd] & \bullet &\bullet\ar@{-}[r] & \bullet \\
\bullet\ar@{-}[rru] & & \bullet\ar@{-}[r] & \bullet \\
\bullet\ar@{-}[rrr] & & & \bullet}\quad\quad
\xymatrix@-5mm{
\bullet\ar@{-}[r] & \bullet\ar@{-}[rd] & \bullet\ar@{-}[r]
 & \bullet \\
\bullet\ar@{-}[rrrd] & & \bullet & \bullet \\
\bullet\ar@{-}[rrru] & & & \bullet}\quad\quad
\xymatrix@-5mm{
\bullet\ar@{-}[r] & \bullet\ar@{-}[r] &\bullet\ar@{-}[r]
 & \bullet \\
\bullet\ar@{-}[rr] & & \bullet\ar@{-}[r] & \bullet \\
\bullet\ar@{-}[rrr] & & & \bullet}
} 
$$
\end{ex}
%
%
\section{Richardson elements}\label{se:rich-theory}
%
%
In this section we describe a method to 
check whether a given nilpotent element of 
the nilradical of a classical Lie algebra 
is a Richardson element. The first statement is 
given in~\cite{bw}. Since we will use 
this result constantly, we repeat its proof. 

\begin{thm}\label{thm:dim-cent-Rich}
Let $\liea{p}\subset\liea{g}$ be a parabolic subalgebra 
of a semi-simple Lie algebra $\liea{g}$, let 
$\liea{p}=\liea{m}\oplus\liea{n}$ where $\liea{m}$ 
is a Levi factor and $\liea{n}$ the corresponding 
nilradical. 
Then $x\in\liea{n}$ is a Richardson element 
for $\liea{p}$ if and only if 
$\dim\liea{g}^x=\dim\liea{m}$. 
\end{thm}

\begin{proof}
Denote the nilradical of the opposite parabolic by 
$\overline{\liea{n}}$ (the opposite parabolic is 
defined as the parabolic subalgebra whose intersection 
with $\liea{p}$ is equal to $\liea{m}$). 
If $x\in\liea{n}$ then 
$\ad(x)\liea{g}=\ad(x)\overline{\liea{n}}+\ad(x)\liea{p}$. Now 
$\ad(x)\liea{p}\subset\liea{n}$ and 
$\dim\ad(x)\overline{\liea{n}}\le\dim\overline{\liea{n}}$. Thus 
\[ 
\dim\ad(x)\liea{g}\le\,2\dim\liea{n}. 
\] 
This implies for $x\in\liea{n}$ that 
$\dim\liea{m}\le \dim\liea{g}^x$ and equality implies that 
$\dim\ad(x)\liea{p}=\dim\liea{n}$. Thus equality implies 
that $x$ is a Richardson element. 

For the other direction, let $x$ be a Richardson 
element for $\liea{p}$. We show that the map 
$\ad(x)$ is injective on $\overline{\liea{n}}$: 
Let $y\in\overline{\liea{n}}$ with $\ad(x)y=0$. 
Then 
\[
0=B(\ad(x)y,\liea{p})=B(y,\ad(x)\liea{p})
=B(y,\liea{n}).
\] 
In particular, $y=0$. 
So $\ad(x)$ is injective on $\overline{\liea{n}}$, 
giving $\dim\ad(x)\overline{\liea{n}}=\dim\liea{n}$. 
Thus 
\begin{eqnarray*}
\dim\overbrace{\ad(x)\liea{p}}^\liea{n}
+
\dim\overbrace{\ad(x)\overline{\liea{n}}}^{\overline{\liea{n}}}
 & = & 2\dim\liea{n} \\ 
 & = & \dim\ad(x)\liea{g} \\ 
 & = & \dim\liea{g}-\dim\liea{g}^x 
\end{eqnarray*}
So 
$\dim\liea{g}^x+\dim\liea{n}=\dim\liea{g}-\dim\liea{n}$
$=\dim\liea{p}=\dim\liea{m}+\dim\liea{n}$, 
i.e. $\dim\liea{m}=\dim\liea{g}^x$. 
\end{proof}

\begin{cor}
Let $\liea{p}=\liea{m}\oplus\liea{n}$ be a parabolic 
subalgebra of a semi-simple Lie algebra. 
Let $X\in\liea{n}$ be a Richardson element.
 
Then $\dim\liea{g}^X\le\dim\liea{g}^Y$ 
for any $Y\in\liea{n}$.
\end{cor}

Theorem~\ref{thm:dim-cent-Rich} gives us a tool to 
decide whether an element of the nilradical of a 
parabolic subalgebra is a Richardson element. 
Namely, we have to calculate its centralizer. 
Centralizers of nilpotent elements of the 
classical Lie algebras can be computed using their  
Jordan canonical form. This 
well-known result is due to Kraft and Procesi, 
cf.~\cite{kp}.

\begin{thm}\label{thm:dim-cent-Jordan} 
Let $(n_1,\dots,n_r)$ be the partition of 
the Jordan canonical form of a nilpotent 
matrix $x$ in the Lie algebra $\liea{g}$, let 
$(m_1,\dots, m_s)$ be the dual partition. 
Then the dimension of the centralizer of $x$ in $\liea{g}$ 
is 
\[
\begin{array}{ll}
\sum\limits_i m_i^2 & \mbox{if $\liea{g}=\liea{gl}_n$}\\
\sum\limits_i \frac{m_i^2}{2}+\frac{1}{2}|\{i\mid n_i\ odd\}|
& \mbox{if }\liea{g}=\liea{sp}_{2n} \\ 
\sum\limits_i \frac{m_i^2}{2}-\frac{1}{2}|\{i\mid n_i\ odd \}| 
 & \mbox{if }\liea{g}=\liea{so}_N
\end{array}
\]
\end{thm}

So it remains to determine the Jordan canonical form 
of a given nilpotent element $x$. It is 
given by the dimensions of the kernels 
of the maps $x^j$, $j\ge 1$: 

\begin{lm}\label{lm:Jordan-form}
Let $x$ be a nilpotent $n\times n$ matrix with $x^{m-1}\neq 0$ 
and $x^m=0$, set $b_j:=\dim\ker x^j$ ($j=1,\dots,m$). 
Define 
\[
a_j:=\left\{ \begin{array}{ll}2b_1-b_2 & j=1 \\ 
2b_j-b_{j-1}-b_{j+1}& j=2,\dots,m-1\\
b_m-b_{m-1} & j=m
\end{array}\right.
\]
Then the Jordan canonical form of $x$ has $a_s$ blocks 
of size $s$ for $s=1,\dots,m$.  
\end{lm}
\begin{cor}\label{cor:part}
With the notation of Lemma~\ref{lm:Jordan-form} above, 
the Jordan canonical form of $x$ is given by the partition 
\[
(1^{a_1},2^{a_2},\dots,(m-1)^{a_{m-1}},m^{a_m}).
\]
\end{cor}

%
\section{The special linear Lie algebra}\label{se:sl-case}
%

We now describe how to obtain a Richardson 
element from a (horizontal) line diagram. 
Recall that a standard parabolic subalgebra of 
$\liea{sl}_n$ is uniquely described by the 
sequence of lengths of the blocks in $\liea{m}$ 
(the standard Levi factor). 
Let $d=(d_1,\dots,d_r)$ be the dimension 
vector of these block lengths. 

We form the horizontal line diagram $L_h(d)$ 
and label its vertices column wise by the numbers 
$1,2,\dots,n$, starting with 
column $1$, labeling top-down. 
This labeled diagram defines a nilpotent element 
as the sum of all elementary matrices $E_{ij}$ 
such that there is a line from $i$ to $j$, 
where $i<j$: 

\[
X(d)=X(L_h(d))
=\sum_{i\mbox{---}j}E_{ij}
\]

\begin{ex}\label{ex:constr}
Let $\liea{p}\subset\liea{sl}_9$ be given by 
the dimension vector $(3,1,2,3)$. We label its  
horizontal line diagram, 
$$
{\small
\xymatrix@-6mm{
1\ar@{-}[r] & 4\ar@{-}[r] & 5\ar@{-}[r] & 7 \\
2\ar@{-}[rr] &  & 6\ar@{-}[r] & 8 \\
3\ar@{-}[rrr] & &  & 9
}},
$$
and obtain $X(d)=$
$E_{1,4}+E_{4,5}+E_{5,7}+E_{2,6}+E_{6,8}+E_{3,9}$, 
an element of the nilradical $\liea{n}$ 
of $\liea{p}$. Using Lemma~\ref{lm:Jordan-form} and 
Corollary~\ref{cor:part} one checks that 
the dimension of the centralizer of $X(d)$ 
is equal to the dimension of the Levi factor. 
Thus $X(d)$ is a Richardson element 
for $\liea{p}$ (by Theorem~\ref{thm:dim-cent-Jordan}). 
\end{ex}
By construction, the matrix $X(d)$ 
is nilpotent for any dimension vector $d$. 
It is in fact an element 
of the nilradical $\liea{n}$ of the parabolic 
subalgebra $\liea{p}=\liea{p}(d)$: 
If $d=(n)$, this is obvious, the 
constructed nilpotent element is the zero matrix. 
If $d=(d_1,d_2)$ then the nonzero 
coefficients of the matrix of $X(d)$ 
are in the rows $1,\dots,d_1$ and columns 
$d_1+1,\dots,d_2$. In other words, they lie in 
the $d_1\times d_2$-block in the upper right 
corner. 
The standard Levi factor consists 
of the blocks $d_1\times d_1$, $d_2\times d_2$ 
on the diagonal. In particular, $X(d_1,d_2)$ 
is a matrix that lies above the Levi factor. This 
generalizes to dimension vectors with more entries.
So we get part (1) of the following Lemma. 
For part (2) we introduce a new notion. 

\begin{defn}
If there exists a sequence of $k$ connected 
lines in a line diagram $L(d)$ that is not 
contained in 
a longer sequence we say that $L(d)$ has 
a {\itshape $k$-chain} or a 
{\itshape chain of length $k$}. A 
{\itshape subchain of length $k$} (or $k$-subchain) 
is a sequence of $k$ connected  
lines in $L(d)$ that maybe contained in a 
longer chain. 
A (sub)chain of length $0$ is a single 
vertex that is not connected to any other vertex. 
\end{defn}
\begin{lm}\label{lm:X-nilrad}
(1) 
The element $X(d)$ is an element 
of the nilradical of $\liea{p}(d)$. 

(2) For $k\ge 1$, the rank of $X(d)^k$ 
is equal to the number of $k$-subchains of lines in 
$L_h(d)$. 
\end{lm}

\begin{proof}[Proof of (2)] 
It is clear that the rank of $X=X(d)$ 
is the number of 
lines in the diagram: to construct $X$, we 
sum over all lines of the diagram. Since these 
lines are disjoint (each vertex $i$ is joint to 
at most one neighbour $j$ with $i<j$) the rows 
and columns of $X$ are linearly independent. 
Therefore the rank of $X$ is equal to the number of 
vertices $i$ such that there is a line from $i$ 
to some $j$ with $i<j$. 

For any $k>0$, the matrix $X^k$ consists of 
linearly independent rows and columns. 
It is clear that an 
entry $(ij)$ of $X\cdot X$ is 
non-zero if and only if there is a line 
$i$---$k$---$j$ in $L_h(d)$: 
$X\cdot X=\sum_{i-k}E_{ik}\sum_{l-j}E_{lj}$ 
where $E_{ik}E_{lj}=\delta_{kl}E_{ij}$. 
Similarly, the rank of $X^k$ is the number 
of vertices 
$i$ such that there exist vertices 
$j_1<j_2<\dots<j_k$ and lines 
$i$---$j_1$---$\,\cdots$---$j_k$ joining 
them,  
i.e. the number of $k$-subchain. 
\end{proof}

It turns out that $X(d)$ is a 
Richardson element for $\liea{p}(d)$, 
as we will show below. This fact follows also 
from the description of Br\"ustle et al. 
in~\cite{bhrr} of $\Delta$-filtered modules 
without self-extension of the Auslander-Reiten 
quiver of type $\lieg{A}_r$ (the number $r$ is 
the number 
of blocks in the standard Levi factor of 
the parabolic subalgebra). 

\begin{thm}\label{thm:lines-rich}
The mapping $d\mapsto X(d)$ 
associates to each dimension vector with $\sum d_i=n$ 
a Richardson 
element for the corresponding parabolic subalgebra 
$\liea{p}=\liea{p}(d)$ of $\liea{sl}_n$. 
\end{thm}

We give here an elementary proof of 
Theorem~\ref{thm:lines-rich} above. 
We will use the ideas of this proof to deal 
with the other classical groups (where we will have 
to use line diagrams that are not horizontal in general). 
The main idea is to use the dimension 
of the centralizer of a Richardson element and the 
partition of the Jordan canonical form of a nilpotent 
element. 

\begin{proof}
Let $d$ be the dimension vector corresponding 
to the parabolic subalgebra $\liea{p}=\liea{p}(d)$. 
Let $X=X(d)$ be the nilpotent element associated 
it (through the horizontal line diagram). 
By Theorem~\ref{thm:dim-cent-Rich} we have to calculate the 
dimension of the centralizer of $X$ and of the Levi 
factor $\liea{m}$ of $\liea{p}$. 
By Theorem~\ref{thm:dim-cent-Jordan}, $\dim\liea{g}^X$ 
is equal to $\sum_i m_i^2-1$ where $(m_1,\dots,m_s)$ 
is the dual partition to the partition of $X$. 
The parts of the dual partition are the entries of 
$d_i$ the dimension vector as is shown in 
Lemma~\ref{lm:diagr-Jordan} below. 
In particular, 
$\dim\liea{l}=\sum_i d_i^2-1=\dim\liea{g}^X$.
\end{proof}

The following result shows how to obtain 
the partition and the dual partition of the Jordan 
canonical form of the nilpotent element 
associated to the dimension vector $d$. 

\begin{lm}\label{lm:diagr-Jordan}
Let $d$ be the dimension vector 
for $\liea{p}\subset\liea{sl}_n$, $X=X(d)$ 
the associated nilpotent element of $\liea{sl}_n$. 
Order the entries $d_1,\dots,d_r$ of the dimension vector 
in decreasing order as 
$D_1,D_2,\dots,D_r$ 
(i.e. such that $D_i\ge D_{i+1}$ for 
all $i$). 
Then the 
Jordan canonical form of $X$ is 
\[
1^{D_1-D_2},2^{D_2-D_3},\dots,(r-1)^{D_{r-1}-D_r},r^{D_r}
\]
and the dual partition is 
\[
D_r,D_{r-1},\dots, D_1.
\]
\end{lm}
In other words, the dual partition for $X(d)$ 
is given by the entries of the dimension vector. 
Furthermore, for every 
$i$-chain in $L_h(d)$ (i.e. for every sequences 
of length $i$, $i\ge 0$, 
that is not contained in a longer sequence) 
the partition has an entry $i$.
\begin{proof}
Let $d=(d_1,\dots,d_r)$ be the dimension 
vector of $\liea{p}$ and $D_1,\dots,D_r$ its permutation 
in decreasing order, $D_i\ge D_{i+1}$. 
To determine the Jordan canonical form of $X=X(d)$ 
we have to compute the rank of the powers $X^s$, $s\ge 1$, 
cf. Lemma~\ref{lm:Jordan-form}. 

Since the nilpotent matrix $X$ is given by the 
horizontal line diagram 
$L_h(d)$, the rank of $X^s$ is 
easy to compute: by Lemma~\ref{lm:X-nilrad} (2), 
the rank of $X^s$ is the number of $s$-subchains. 
In particular, 
$\rk X=n-D_1$ and $\rk X^2=n-D_1-D_2$, 
$\rk X^3=n-D_1-D_2-D_3$, etc. 
This gives 
\[
b_s:=\dim\ker X^s=D_1+\dots+D_s \ \mbox{for} \ s=1,\dots,r.
\] 
And so, by Lemma~\ref{lm:Jordan-form}, we obtain 
$a_1=D_1-D_2$, $a_2=D_2-D_3$, $\dots,a_r=D_r$ 
proving the first statement. 
The statement about the dual partition (i.e. 
the partition given by the lengths of the columns 
of the partition) follows then immediately. 
\end{proof}

%
%
\section{Richardson elements for the other 
classical Lie algebras}\label{se:BCD-type}
%
%
In this section we will introduce generalized 
line diagrams to deal with the symplectic and 
orthogonal Lie algebras. Having introduced 
them, we show that they correspond to 
Richardson elements for the parabolic subalgebra 
in question. Then we discuss some properties 
and describe the dual of the partition of 
a nilpotent element given by such a generalized 
line diagram. 
Furthermore, we describe the support of the 
constructed $X(d)$ and relate it to the 
Bala-Carter label of the $G$-orbit through 
$X(d)$ where $G$ is the adjoint group of $\liea{g}$. 

To define the orthogonal Lie 
algebras, we use the skew diagonal matrix $J_n$ 
with ones on the skew diagonal and zeroes 
else. 
The symplectic Lie algebras $\liea{sp}_{2n}$ 
are defined using 
${\small\begin{bmatrix} 0 & J_n \\ -J_n & 0\end{bmatrix}}$. 
(For details we refer the reader to~\cite{gw}.) 
So 
$\liea{so}_n$ consists of the $n\times n$-matrices 
that are skew-symmetric around the skew-diagonal 
and $\liea{sp}_{2n}$ is the set of $2n\times 2n$-matrices 
of the form 
\[\begin{bmatrix}A & B\\ C&A^*\end{bmatrix}\]
where $A^*$ is the the negative of the skew transpose 
of $A$. 
Thus in the case of the symplectic and orthogonal 
Lie algebras, the block sizes of the standard Levi 
factor form a palindromic sequence. 

If there is an even number of blocks in the 
Levi factor, the dimension vector is of the form 
$(d_1,\dots,d_r,d_r,\dots,d_1)$. 
We will refer to this situation as type~(a). 
If there is an odd number of blocks in the 
Levi factor, type (b), the dimension vector 
is $(d_1,\dots,d_r,d_{r+1},d_r,\dots,d_1)$. 

By the (skew) symmetry around the skew diagonal,  
the entries below the skew diagonal 
of the matrices $X(d)$ are 
determined by the entries above the skew 
diagonal. In terms of line diagrams: 
For $\liea{sp}_N$ and $\liea{so}_N$ there is 
a line $(N-j+1)$---$(N-i+1)$ whenever 
there is a line $i$---$j$. 
We will call the line $(N-j+1)$---$(N-i+1)$ 
the {\itshape counterpart} of $i$---$j$ and 
will sometimes denote counterparts by dotted lines.
In particular, it suffices to describe 
the lines attached to the left to vertices of the 
first $r$ 
columns for both types (a) and (b). 

The (skew)-symmetry will give constraints on the 
diagram - there will also appear negative 
entries. For the moment, let us assume 
that $L(d)$ is a diagram defining an element 
of the nilradical of the parabolic subalgebra 
in question. Then part (2) of 
Lemma~\ref{lm:X-nilrad} still holds.
\begin{lm}\label{lm:chains-rank}
If $X(d)$ is defined by $L(d)$
then the rank of the 
map $X(d)^k$ is the number of $k$-subchains 
of lines in the diagram. 
\end{lm}
This uses the same argument as 
Lemma~\ref{lm:X-nilrad} since by construction, 
$X(d)$ only has linearly independent 
rows and columns and the product 
$X(d)^2$ only has nonzero entries $E_{il}$ 
if $X(d)$ has an entry $E_{ij}$ and an 
entry $E_{jl}$ for some $j$. 

The following remark allows us to simplify the 
shapes of the diagrams we are considering. 
If $d=(d_1,\dots,d_r)$ is an 
$r$-tuple in $\NN^r$, and $\sigma\in S_r$ (where 
$S_r$ is the permutation group on $r$ letters) 
we define $d_{\sigma}$ 
as $(d_{\sigma 1},d_{\sigma 2},\dots,d_{\sigma r})$. 
By abuse of notation, for 
$d=(d_1,\dots,d_r,d_r,\dots,d_1)$ in $\NN^{2r}$,  
we write 
$d_{\sigma}=(d_{\sigma 1},\dots,d_{\sigma r},
d_{\sigma r},\dots,d_{\sigma 1})$ 
and for $d=(d_1,\dots,d_r,d_{r+1},d_r,\dots,d_1)$ 
in $\NN^{2r+1}$, we define 
$d_{\sigma}$ to be the $2r+1$-tuple 
$(d_{\sigma 1},\dots,d_{\sigma r},
d_{r+1},d_{\sigma r},\dots,d_{\sigma 1})$. 
It will be clear from the context which tuple 
we are referring to. 
\begin{re}\label{re:permutations}
For $d=(d_1,\dots,d_r)$ the diagrams 
$L_h(d)$ and $L_h(d_{\sigma})$ 
have the same chains of lines for any $\sigma\in S_r$. 
In other words: for any $k\ge 1$, the number of 
chains of lines of length $k$ in $L_h(d)$ is the 
same as the number of lines of length $k$ 
in $L_h(d_{\sigma})$.  
As an illustration, consider the permutation $1243$ of 
$d=(3,1,2,3)$: 
$$
{\small
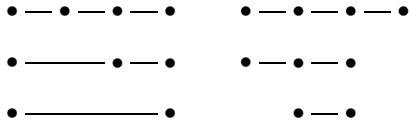
\xymatrix@-5mm{
\bullet\ar@{-}[r] & \bullet\ar@{-}[r] &\bullet\ar@{-}[r]
 & \bullet \\
\bullet\ar@{-}[rr] & & \bullet\ar@{-}[r] & \bullet \\
\bullet\ar@{-}[rrr] & & & \bullet}\quad\quad
\xymatrix@-5mm{
\bullet\ar@{-}[r] & \bullet\ar@{-}[r] &\bullet\ar@{-}[r]
 & \bullet \\
\bullet\ar@{-}[r] & \bullet\ar@{-}[r] &\bullet  &  \\
 & \bullet\ar@{-}[r] & \bullet &}
} 
$$
Similarly, for $f=(f_1,\dots,f_r,f_r,\dots,f_1)$ 
resp. for $g=(g_1,\dots,g_r,g_{r+1},g_r,\dots,g_1)$, 
if $L(f)$ and $L(g)$ are line diagrams for $\liea{sp}_{2n}$ 
or $\liea{so}_N$ then for any 
$\sigma\in S_r$, the diagrams 
$L(f_{\sigma})$ resp. $L(g_{\sigma})$ are 
also diagrams for the corresponding Lie algebras 
and have the same exactly the same chains as 
$L(f)$ resp. as $L(g)$. 
\end{re}
We have an immediate consequence of 
Remark~\ref{re:permutations} and of 
Lemma~\ref{lm:chains-rank}: 

\begin{cor}\label{cor:reordering}
Let $d=(d_1,\dots,d_r,d_r,\dots,d_1)$ or 
$d=(d_1,\dots,d_r,d_{r+1},d_r,\dots,d_1)$ 
be the dimension vector of a parabolic 
subalgebra of a symplectic or orthogonal 
Lie algebra and $X(d)$ be given by the 
appropriate line diagram. 
In calculating the rank of $X(d)^k$ we 
can assume that $d_1\le\dots\le d_r$. 
\end{cor}

We will make frequent use of this 
property. 
Now we will finally be able to construct 
diagrams for the other classical cases. 
We have already mentioned that the horizontal 
line diagrams do not produce Richardson 
elements. 
One reason is that the counterpart of 
a line $i$---$j$ is not always horizontal. 
The other reason is that we have to introduce 
negative signs for the symplectic and 
orthogonal cases when we associate a 
nilpotent matrix to a diagram:  
If $\liea{g}=\liea{sp}_{2n}$, in the 
definition of $X(d)$ 
we subtract $E_{ij}$ whenever there is a 
line $i$---$j$ with $n<i<j$. 
If $\liea{g}=\liea{so}_N$ we subtract 
$E_{ij}$ whenever there is a line $i$---$j$ 
with $i+j>N$. 

\begin{ex}\label{ex:sp-non-horizontal}
Let $(1,2,2,1)$ be the dimension vector 
of a parabolic subalgebra of $\liea{sp}_6$. 
Then the following three line diagrams 
determine elements of the nilradical 
of $\liea{p}$: 
$${\small
\xymatrix@-6mm{
1\ar@{-}[r] & 2\ar@{-}[r] & 4 & 6\\
        & 3\ar@{.}[r] & 5\ar@{.}[ur] \\
}\quad\quad
\xymatrix@-6mm{
1\ar@{-}[r] & 2\ar@{-}[rd] & 4 & 6\\
        & 3\ar@{-}[ru] & 5\ar@{.}[ur] \\
}\quad\quad
\xymatrix@-6mm{
1\ar@{-}[r] & 2\ar@{-}[r] & 5\ar@{.}[r] & 6\\
        & 3\ar@{-}[r] & 4 \\
}}
$$ 
The last diagram is just a reordering 
of the second. 
The nilpotent elements are 
$X_1=E_{12}+E_{24}+E_{35}-E_{56}$ 
resp. $X_2=E_{12}+E_{25}+E_{34}-E_{56}$. 
By calculating the Jordan canonical 
forms for these elements one 
checks that only the nilpotent element 
$X_2$ is a 
Richardson element.
\end{ex}

This example and the discussion above illustrate 
that for the symplectic and orthogonal Lie 
algebras, we will use:

(i) non-horizontal lines,

(ii) labeling top-bottom {\bf and} bottom-top, 

(iii) negative signs, too.

\noindent
Before we start defining these line diagrams
we introduce a new notion. 
\begin{defn}
Let $\liea{p}$ be the standard parabolic 
subalgebra of a symplectic or orthogonal 
Lie algebra $\liea{g}$. 
We say that $\liea{p}$ is {\itshape simple} 
if $\liea{p}\subset\liea{g}$ 
is of one of the following forms: 
\begin{enumerate}
\item A parabolic subalgebra of $\liea{sp}_{2n}$ 
with an even number of 
blocks in the standard Levi factor.
\item A parabolic subalgebra of $\liea{so}_{2n}$ 
with an even number of 
blocks in the standard Levi factor such that 
odd block lengths appear exactly twice. 
\item
A parabolic subalgebra of $\liea{sp}_{2n}$ 
with an odd number of 
blocks in the Levi factor and such that each 
odd $d_i$ that is smaller than $d_{r+1}$ 
appears exactly twice. 
\item A parabolic subalgebra of $\liea{so}_N$ 
with an odd number of 
blocks in the Levi factor such that 
either all $d_i$ are odd 
or there is an index $k\le r$ such 
that all $d_i$ with $i\le k$ are even, 
$d_j$ odd for $j>k$ and the even 
$d_i$ are smaller than $d_{k+1},\dots,d_r$. 
Furthermore, the even block lengths that 
are larger than $d_{r+1}$ appear only
once among $d_1,\dots,d_k$. 
\end{enumerate}
\end{defn}

\begin{defn}[Type (a)]
Let $\liea{p}$ be a simple parabolic subalgebra 
of $\liea{sp}_{2n}$ or $\liea{so}_{2n}$, 
given by the dimension vector 
$d=(d_1,\dots,d_r,d_r,\dots,d_1)$. 
Then we define the 
{\itshape line diagram} $L_{even}(d)$ 
{\itshape associated to} $d$ 
(and $\liea{g}$) 
as follows. 

\begin{enumerate}
\item
Draw $2n$ vertices in $2r$ columns of length 
$d_1,\dots$, top-adjusted. 
Label the first $r$ columns with 
the numbers $1,\dots, n$, top--bottom. 
Label the second $r$ columns with the numbers 
$n+1,\dots, 2n$, bottom--top. 
\item
Join the first $r$ columns with horizontal lines 
as for $\liea{sl}_n$. Draw the counterparts of 
these lines in the second $r$ columns. 
\item[(3) (i)]
If $\liea{g}=\liea{sp}_{2n}$, add 
the lines $k$---$(2n-k+1)$. 
\item[(3) (ii)]
If $\liea{g}=\liea{so}_{2n}$, one adds 
the lines $(2l-1)$---$(2n-2l+1)$ and their 
counterparts $2l$---$(2n-2l+2)$ if $n$ 
is even. If $n$ is odd, the lines 
$2l$---$(2n-2l)$ and their counterparts 
$(2l+1)$---$(2n-2l+1)$.
\end{enumerate}
\end{defn}

\begin{defn}[Type (b)]
Let $\liea{p}$ be a simple parabolic subalgebra 
of $\liea{sp}_{2n}$ or of $\liea{so}_N$, 
given by the dimension vector 
$d$ 
$=(d_1,\dots,d_r,d_{r+1},d_r,\dots,d_1)$. 
Then we define the 
{\itshape line diagram} $L_{odd}(d)$ 
{\itshape associated to} 
$d$ (and $\liea{g}$) 
as follows. 

\begin{enumerate}
\item
Draw $2r+1$ columns of length $d_1,\dots$, 
top-adjusted. Label them with the numbers 
$1,\dots$ in 
increasing order, top--bottom in each column. 
\item[(2) (i)]
For $\liea{sp}_{2n}$: \\
If $\min_i\{d_i\}\ge 2$, draw a horizontal 
of lines in the first row and all their 
counterparts, forming a sequence joining the lowest 
vertices of each column. Repeat this procedure 
as long as the columns of the remaining 
vertices are all at least of length two. 
\item[(2) (ii)]
For $\liea{so}_N$: \\
If $d_1$ is odd, go to step (3) (ii). 
If $d_1$ is even, do as in (2) (i), 
drawing lines in the first row and their 
counterparts joining the lowest vertices. 
Repeat until either the first of the remaining 
columns has odd length or there are no 
vertices left to be joined. 
Continue as in (3) (ii). 
\item[(3) (i)]
For $\liea{sp}_{2n}$: \\
For the remaining vertices: draw horizontal 
lines following the top-most remaining 
vertices and simultaneously their 
counterparts (the lowest remaining 
vertices). 
\item[(3) (ii)]
For $\liea{so}_N$: \\
All columns 
have odd length. Connect the central 
entries of each column. The remaining 
column lengths are all even, the are 
joined as in (2) (ii). 
\end{enumerate}
\end{defn}

\begin{thm}\label{thm:line-richardson}
Let $d$ be the dimension vector 
for a simple parabolic subalgebra of $\liea{sp}_{2n}$ 
or $\liea{so}_N$. 
Then the associated diagram $L_{even}(d)$ 
resp. 
$L_{odd}(d)$ determines a Richardson 
element for $\liea{p}(d)$ by setting 
\[
\begin{array}{ccll}
X(d) & = & \sum_{i\mbox{---}j,\  i\le n}E_{ij} 
 - \sum_{i\mbox{---}j,\ i>n}E_{ij}
 & \mbox{for}\  \liea{sp}_{2n}\\ 
X(d) & = & \sum_{i\mbox{---}j,\ i+j<N}E_{ij} 
 - \sum_{i\mbox{---}j,\ i+j>N}E_{ij} 
 & \mbox{for} \ \liea{so}_N
\end{array}
\]
where the sums are over all lines in the diagram.
\end{thm}
We first include some 
immediate consequences of this result. 
After that we add an observation about 
the (dual of the) partition corresponding to $X(d)$ 
and then we are ready to prove 
Theorem~\ref{thm:line-richardson}. 

Theorem~\ref{thm:line-richardson} 
enables us to 
determine the minimal 
$k$ such that the Richardson element $X(d)$ 
lies in the graded parts 
$\liea{g}_1\oplus\dots\oplus\liea{g}_k$. 
To do so we introduce $s(d)$ 
as the maximal number of 
entries $d_i,\dots,d_{i+s}$ 
of $d$ that are surrounded by larger 
entries $d_{i-1}$ and $d_{i+s+1}$. 
More precisely, if $d=(d_1,\dots,d_r,d_r,\dots,d_1)$ 
or $d=(d_1,\dots,d_r,d_{r+1},\dots,d_1)$ 
is the dimension vector, we rewrite $d$ 
as a vector with increasing indices, 

$(c_1\dots,c_r,c_{r+1},c_{r+2},\dots,c_{2r})$ 
resp. 

$(c_1\dots,c_r,c_{r+1},c_{r+2},\dots,c_{2r+1})$ 
and define 

$s(d)
:=1+\max_i\{\text{there are}\ c_{j+1},\dots, c_{j+i}\mid 
c_j>c_{j+l}<c_{j+i+1}\text{ for all }0\le l\le i\}$.

\begin{cor}\label{cor:bound-grade}
Let $\liea{p}(d)$ be a simple parabolic 
subalgebra of the orthogonal or symplectic 
Lie algebras. 
Then the element $X(d)$ belongs to 
$\liea{g}_1\oplus\dots\oplus\liea{g}_{s(d)}$. 
The same holds for parabolic subalgebras 
of $\liea{sl}_n$. 
\end{cor}
This follows from the fact that $E_{ij}$ 
with $i$ from column $k$ of the line diagram 
and $j$ from column $k+s$
is an entry of the graded part $\liea{g}_s$. 
If, e.g., we have $c_1>c_j<c_{s+1}$ for $j=2,\dots,s$ 
then there is a line joining 
columns one and $s+1$. So $X(d)$ 
has an entry in $\liea{g}_{s}$. 
\begin{cor}
For $\liea{sl}_n$, $s(d)$ is equal to one 
if and only if the dimension vector satisfies 
$d_1\le\dots\le d_t\ge\dots\ge d_r$ for some 
$1\le t\le r$. 
\end{cor}

This well-known result has been observed by 
Lynch~\cite{l}, Elashvili and Kac~\cite{ek}, 
Goodwin and R\"ohrle~\cite{gr}, 
and in our joint work with Wallach~\cite{bw}. 

The next lemma shows how to obtain the dual of 
the partition of $X(d)$ if $X(d)$ is given by 
the appropriate line diagram for $d$. 

\begin{lm}\label{lm:dual-part}
If $\liea{p}(d)$ is a simple parabolic subalgebra of a 
symplectic or orthogonal Lie algebra let 
$X=X(d)$ be given by the appropriate line diagram 
$L_{even}(d)$ or $L_{odd}(d)$. The 
dual of the partition of $X$ has the form 
\[
\begin{array}{llll}
 & \text{Dual of the partition of $X$} & \liea{g}
 & \text{Type of $\liea{p}$} \\
 & \\
(i) & d_1,d_1,\dots,d_r,d_r & \liea{sp}_{2n} & (a) \\
 & & \\ 
(ii) & d_{r+1}\cup 
\left(\bigcup_{d_i\notin D_o}d_i,d_i\right)
\cup\left(\bigcup_{d_i\in D_o} d_i-1,d_i+1\right)
 & \liea{sp}_{2n} & (b) \\
 & & \\
(iii) & \left(\bigcup_{d_i\text{even}} d_i,d_i\right)
\cup\left(\bigcup_{d_i\text{odd}}d_i-1,d_i+1\right)
 & \liea{so}_{2n} & (a) \\
 & & \\
(iv) & d_{r+1}\cup\left(
\bigcup_{d_i\notin D^e}
d_i,d_i\right)
\cup\left(\bigcup_{d_i\in D^e}d_i-1,d_i+1\right)
 & \liea{so}_{2n+1} & (b) \\
 &  \\
(v) & d_{r+1}\cup\left(
\bigcup_{d_i\notin D^o}
d_i,d_i\right)
\cup\left(\bigcup_{d_i\in D^o}d_i-1,d_i+1\right)
 & \liea{so}_{2n} & (b) 
\end{array}
\]
where 
$D_o:=\{d_i\text{ odd}\mid d_i<d_{r+1}\}$, 
$D^o:=\{d_i\text{ odd}\mid d_i>d_{r+1}\}$ 
and 
$D^e:=\{d_i\text{ even}\mid d_i>d_{r+1}\}$ 
are subsets of $\{d_1,\dots,d_r\}$. 

In particular, if $D_o$, $D^e$ or $D^o$ are 
empty, the partition in the corresponding 
case (ii), (iv) or (v) has the same parts 
as the dimension vector. The same is true 
for (iii), if there are no odd $d_i$. 
\end{lm}
The proof consists mainly in 
counting lines and (sub)chains of lines of 
the corresponding diagrams. Therefore we 
postpone it and include it in the appendix.
We are now ready to prove Theorem~\ref{thm:line-richardson} 
with the use of Theorem~\ref{thm:dim-cent-Jordan} 
and of Lemma~\ref{lm:dual-part}.
\begin{proof}[Proof of Theorem~\ref{thm:line-richardson}]
We consider the case $\liea{g}=\liea{sp}_{2n}$. 
For the parabolic subalgebras of an orthogonal Lie algebra, 
the claim follows using the same methods. 
The idea is to use the 
dimension of the centralizer 
of $X(d)$ and compare it to the dimension of 
the Levi factor. To calculate the dimension of 
the centralizer, we use the formulae of 
Theorem~\ref{thm:dim-cent-Jordan}, i.e. we use the 
dual of the partition of $X=X(d)$ as described 
in Lemma~\ref{lm:dual-part} and the number of 
odd parts in the partition of $X$. 

\noindent
$\liea{sp}_{2n}$, type (a): \\
By Lemma~\ref{lm:dual-part}
the dual partition of the nilpotent 
element $X=X(d)$ has as parts the entries of $d$. 
Since they all appear in pairs, the partition 
of the orbit has no odd entries. So 
by the formula of 
Theorem~\ref{thm:dim-cent-Jordan} we 
obtain 
$\dim\liea{g}^X=\frac{1}{2}(2d_1^2+\dots+2d_r^2)$, 
the same as the dimension of the Levi 
factor. In particular, $X$ is a Richardson element 
for the parabolic subalgebra $\liea{p}(d)$ of $\liea{sp}_{2n}$. 

\noindent
$\liea{sp}_{2n}$, type (b): \\
As in Lemma~\ref{lm:dual-part} let 
$D_o\subset\{d_1,\dots,d_r\}$ be the possibly 
empty set of the odd $d_i$ that are smaller than 
$d_{r+1}$. Then the dual partition has the parts 
\[
\{d_i,d_i \mid i<r,\ d_i\notin D_o\}\cup\{d_{r+1}\}
\cup\{d_{i+1},d_{i-1}\mid d_i\in D_o\}.
\]
The $d_i$ that are not in $D_o$ come in pairs and do not 
contribute to odd parts in the partition of 
$X=X(d)$. In particular, the number of odd parts only 
depends on $d_{r+1}$ and on the entries of $D_o$. 
We write the 
elements of $D_o$ in decreasing order as 
$\tilde{d}_1,\dots,\tilde{d}_s$ (where $s=|D_o|$). 
By assumption (the parabolic subalgebra is simple)
these odd entries are all different, 
$\tilde{d}_1>\tilde{d}_2>\dots>\tilde{d}_s$. 
Then the number of odd parts of the partition 
of $X$ is the same as the number of odd parts 
of the dual of the partition
\[
\tilde{P}:\quad d_{r+1},\tilde{d}_1+1,\tilde{d}_1-1,\dots,
\tilde{d}_s+1,\tilde{d}_s-1.
\]
This has $d_{r+1}-(\tilde{d}_1+1)$ ones, 
$(\tilde{d}_1+1)-(\tilde{d}_1-1)$ twos, 
$(\tilde{d}_1-1)-(\tilde{d}_2+1)$ threes, and so on. 
So the number of odd parts in the dual of 
$\tilde{P}$ is 
\[ 
[d_{r+1}-(\tilde{d_1}+1)]+[(\tilde{d}_1-1)-(\tilde{d}_2+1)] 
 + \dots
 + [(\tilde{d}_{s-1}-1)-(\tilde{d}_s+1)]+\tilde{d}_s-1 \\ 
 = d_{r+1}-2s.
\]
Thus the dimension of the centralizer of $X$ is 
\begin{align*}
\frac{1}{2} 
& \left[
\left(\sum_{\substack{i<r+1\\ d_i\notin D_o}}2d_i^2\right)
+ d_{r+1}^2 + 
\left(\sum_{d_i\in D_o}(d_i-1)^2+(d_i+1)^2\right) 
 + d_{r+1}-2s\right] \\ 
&  = \sum_{i\le r}d_i^2
+\binom{d_{r+1}+1}{2}  = \dim\liea{m}.
\end{align*}
\end{proof}

\subsection{Bala Carter labels for Richardson 
orbits}

The support of the nilpotent element of a 
simple line diagram is by construction a 
simple system of root. Namely, for any $d$, 
the corresponding $X(d)$ has 
at most one non-zero element in each row 
and each column. 
One can check that none of the corresponding 
positive roots subtract from each other. 

In other words, the support 
$\supp(X)$ forms a simple system of roots.  
\begin{re}
The converse statement is not true. 
There are Richardson elements whose support 
form a simple system of roots but where 
there is no simple line diagram defining 
a Richardson element.
A family of examples are the Borel subalgebras 
of $\liea{so}_{2n}$ or more general, 
parabolic subalgebras of $\liea{so}_{2n}$ 
where 
$\alpha_n$ and $\alpha_{n-1}$ are both not 
roots of the Levi factor
\end{re}

If $X$ is a nilpotent element of $\liea{g}$ 
we denote 
the $G$-orbit through $X$ by $\mathcal{O}_X$ 
(where $G$ is the adjoint group of $\liea{g}$). 

\begin{cor}
Let $\liea{p}(d)$ be a 
parabolic subalgebra of $\liea{sl}_n$. 
Define $X(d)$ by the 
line diagram $L_h(d)$ 
or a simple parabolic subalgebra of 
(b)-type for $\liea{sp}_{2n}$, $\liea{so}_N$ 
Then the group spanned by $\supp X(d)$ is equal to 
the Bala-Carter label of the $G$-orbit 
$\mathcal{O}_{X(d)}$. 
\end{cor}

\begin{proof}
This follows from the characterization 
of the type (i.e. the Bala-Carter label) 
of $\mathcal{O}_X$ given by 
Panyushev in Section 3 of~\cite{pan}. 

For simplicity we assume $d_1\le \dots\le d_r$. 
Note that in any case, the partition of 
$X(d)$ is given by the chains in the 
line diagram. The partition of $X(d)$ has 
entry $i$ for every chain of length $i+1$. 

If $\alpha$ given by $E_{ij}$ and $\beta$ given 
by $E_{kl}$ are roots of 
$\supp X(d)$ then they add to a root of 
$\liea{sl}_n$ if and only if there is a
line connecting them. 
Thus in the case of the special linear Lie algebra a chain 
of length $i+1$ corresponds to a factor 
$\lieg{A}_i$ in $\supp X(d)$. 
Similarly, for $\liea{sp}_{2n}$ and $\liea{so}_N$, 
a chain of 
length $i+1$ together with its counterpart give 
a factor $\lieg{A}_i$. 
Finally, the possibly remaining 
single chain of length $2j+1$ (passing 
through the central vertex of column $r+1$) 
in the case of $\liea{so}_{2n+1}$ 
gives a factor $\lieg{B}_j$. 
Then the claim follows with~\cite{pan} where 
Panyushev describes the type of a nilpotent orbit 
in terms of its partition. 
\end{proof}

%
\section{Branched diagrams}\label{se:branched}
%

The diagrams we have introduced had at 
most one line to the left and at most one 
line to the right of a vertex. 
We call such a diagram a 
{\itshape simple line diagram}. 
In the case of simple parabolic subalgebras, 
we can always choose 
a simple line diagram to define a 
Richardson element. 
However, there are parabolic subalgebras 
where no simple diagram gives rise to 
a Richardson elements. 
After giving an example we characterize 
the parabolic subalgebras for which 
there exists a simple line diagram giving 
a Richardson element. 
Then we discuss the case of the symplectic 
Lie algebras. We introduce a branched 
diagram and obtain a 
Richardson elements for the parabolic 
subalgebra in question. 

\begin{ex}
1) Consider the parabolic subalgebra of 
$\liea{so}_{2n}$ given by the dimension 
vector $(n,n)$ where $n$ is odd. 
The element $X=X(n,n)$ 
given by the diagram $L_{even}(n,n)$ 
has rank $n-1$ and so the kernel of the 
map $X^k$ has dimension $n+1$ or $2n$ 
for $k=1,2$ resp. The partition of $X$ 
is then $1^2,2^{n-1}$, its dual 
is $n-1,n+1$. The centralizer of 
$X$ has dimension $2n^2+1-1$ and the 
Levi factor of this parabolic subalgebra 
has dimension $n^2$. So $X$ is a Richardson 
element.

2) Let $\liea{p}\subset\liea{so}_{4d}$ 
be given by $(d,d,d,d)$ where $d$ is odd. 
Note that the skew-symmetry of the orthogonal 
Lie algebra allows at most $d-1$ lines between 
the two central columns. 
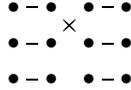
\[
{\small
\xymatrix@-7mm{
\bullet\ar@{-}[r] & \bullet\ar@{-}[rd] & \bullet\ar@{-}[r] 
& \bullet \\
\bullet\ar@{-}[r] & \bullet\ar@{-}[ru] & \bullet\ar@{-}[r]
 & \bullet \\ 
\bullet\ar@{-}[r] &\bullet & \bullet\ar@{-}[r] & \bullet
}}
\]
The line diagram $L_{even}(d,d,d,d)$ 
has $2d+d-1$ lines, 
$2(d-1)$ two-subchains and 
$d-1$ three-chains. Calculating the 
dimensions of the kernel of the 
map $X^k$ (where $X=X(d,d,d,d)$) yields the partition 
$2^2,4^{d-1}$. Its dual is 
$(d-1)^2,(d+1)^2$, hence the centralizer 
of $X$ has dimension $2d^2+2$ while 
the Levi factor has dimension $2d^2$. 
\end{ex}

\begin{thm}
Let $\liea{g}$ be a simple Lie algebra. The 
parabolic subalgebras $\liea{p}$ of 
$\liea{g}$ for which 
there exists a simple line diagram that defines 
a Richardson element for $\liea{p}$ are: 

The parabolic subalgebras of $\liea{sl}_n$ and 
the simple parabolic subalgebras of the symplectic 
and orthogonal Lie algebras.
\end{thm}

\begin{proof}
By Theorems~\ref{thm:lines-rich} 
and~\ref{thm:line-richardson} 
there is always a simple line diagram giving 
a Richardson element in these cases. 
It remains to show that these are the only ones. 
By Corollary~\ref{cor:reordering} 
we can assume w.l.o.g. that $d_1\le\dots\le d_r$. 
Then it turns out that if there is an even number 
of blocks for $\liea{so}_{2n}$ or 
if $d_r\le d_{r+1}$ for $\liea{sp}_{2n}$ 
the problem is translated to the problem of finding 
a Richardson element in the first graded 
part $\liea{g}_1$ of $\liea{g}$ because 
of the following observation: 
Since 
$d_1\le\dots\le d_r=d_r\ge\dots\ge d_1$, 
or $d_1\le\dots\le d_r\le d_{r+1}\ge d_r\ge\dots\ge d_1$ 
all lines are connecting neighbored 
columns. But lines connecting neighbored 
columns correspond to entries $E_{i,j}$ 
of the first super diagonal of the parabolic 
subalgebra, i.e. to entries of $\liea{g}_1$. 
Then the claim follows from the classification 
of parabolic subalgebras with a Richardson 
element in $\liea{g}_1$ for 
type (a) of $\liea{so}_{2n}$ 
and if $d_r\le d_{r+1}$ for type (b) 
parabolic subalgebras of the symplectic 
Lie algebra. 
In both cases there exists a Richardson 
element in $\liea{g}_1$ if and only if 
each odd block length $d_i$ only appears 
once among $d_1,\dots,d_r$, cf.~\cite{bw}. 
If there is no Richardson element in 
$\liea{g}_1$ then in particular no simple line 
diagram can give a Richardson element. 
It remains to deal with (b)-types for $\liea{so}_N$ 
and (b)-types for $\liea{sp}_{2n}$ 
where $d_{r+1}$ is not maximal. 
Both are straightforward but rather lengthy calculation 
that we omit here.
\end{proof}

By way of illustration we include examples of 
branched diagrams for non-simple parabolic 
subalgebras of $\liea{sp}_{2n}$ and of $\liea{so}_N$ 
in the appendix. 
In general, it is not clear how branched diagrams 
should be defined uniformely for the 
symplectic and orthogonal Lie algebras. 
It is clear from the description of simple 
parabolic subalgebras of $\liea{so}_N$ that 
this case is more intricate. We assume that 
Richardson elements can be obtained by 
adding lines to the corresponding simple 
line diagrams: 

\begin{conj}
For the (b)-type of $\liea{sp}_{2n}$ 
the appropriate diagram defining a Richardson 
element is obtained from $L_{odd}(d)$ 
by adding a branching 
for every repetition $d_i=d_{i+1}=\dots=d_{i+s}$ 
of odd entries smaller than $d_{r+1}$. 
\end{conj}

We conclude this section with a remark on 
the bound $s(d)$ introduced in 
Section~\ref{se:BCD-type}. If there is 
no simple line diagram defining a Richardson 
element, we can still define 
$s(d)$ to be the maximal number of a sequence 
of entries of $d$ that are surrounded by 
two larger entries. But this will now 
only be a lower bound, the Richardson 
element defined by a branched diagram 
does not necessarily lie in 
$\liea{g}_1\oplus\dots\oplus\liea{g}_{s(d)}$, 
cf. Examples~\ref{ex:branched-sp},~\ref{ex:branched-sp22}, 
and~\ref{ex:branched-so}. 

%
%
%
\section*{Appendix}
%
%
We discuss some examples of branched line diagrams 
for $\liea{sp}_{2n}$ and for $\liea{so}_N$  
to illustrate Section~\ref{se:branched}. 
Recall that the parabolic subalgebras of 
type (b) of $\liea{sp}_{2n}$ are simple if 
and only if every odd $d_i<d_{r+1}$ only appears 
once among $d_1,\dots,d_r$. In particular, the 
smallest example of $\liea{sp}_{2n}$ where there is no 
simple line exists for $n=3$. 

\begin{ex}\label{ex:branched-sp}
Let $\liea{p}$ be the 
parabolic subalgebra of $\liea{sp}_6$ with 
dimension vector $(1,1,2,1,1)$. Consider 
the diagrams 
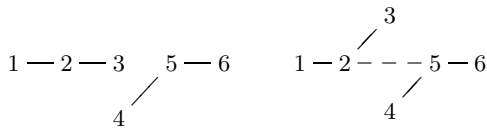
\[
{\small
\xymatrix@-5mm{ & & \ &  \\
1\ar@{-}[r] & 2\ar@{-}[r] & 3 & 5\ar@{-}[r] & 6 \\ 
 & & 4\ar@{-}[ru]
}\quad\quad 
\xymatrix@-6mm{ & & 3 \\
1\ar@{-}[r] & 2\ar@{-}[ru]\ar@{--}[rr] & 
 & 5\ar@{-}[r] & 6 \\ 
 & & 4\ar@{-}[ru]
}}
\]
The diagram to the left is a line diagram as in 
Section~\ref{se:BCD-type}. The corresponding 
nilpotent element has a centralizer of dimension 
$7$. 
However, the Levi factor is five dimensional. 
In the second diagram, there is one extra 
line, connecting the vertices $2$ and $5$. The defined 
matrix $X=E_{12}+E_{23}+E_{25}-E_{45}-E_{56}$ 
has a five dimensional centralizer 
as needed. 
\end{ex}

\begin{ex}\label{ex:branched-sp22}
The following branched line diagram for the 
parabolic subalgebra of 
$\liea{sp}_{22}$ with dimension vector 
$d=(1,1,1,3,3,4,3,3,1,1,1)$ gives a Richardson 
element for $\liea{p}(d)$ 
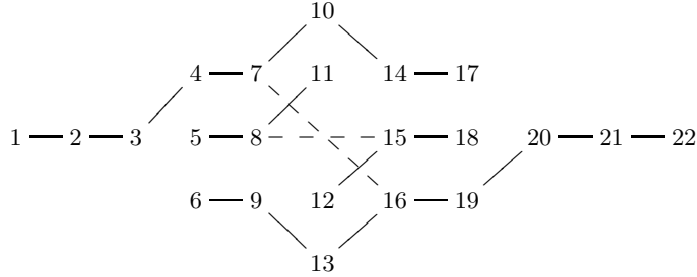
\[
{\small
\xymatrix@-4mm{ & & & & & 10\ar@{-}[rd] \\
 & & & 4\ar@{-}[r] & 7\ar@{-}[ru]\ar@{--}[rrdd]
 & 11 & 14\ar@{-}[r] & 17 \\
1\ar@{-}[r] & 2\ar@{-}[r] & 3\ar@{-}[ru] & 5\ar@{-}[r]
 & 8\ar@{-}[ru]\ar@{--}[rr]
 &  & 15\ar@{-}[r]
 & 18 & 20\ar@{-}[r] & 21\ar@{-}[r] & 22 \\
 & & & 6\ar@{-}[r] & 9\ar@{-}[rd] & 12\ar@{-}[ru]
 & 16\ar@{-}[r] & 19\ar@{-}[ru] & \\
 & & &  &  & 13\ar@{-}[ru] 
}}
\]
The Levi factor and the centralizer of the 
constructed $X$ have dimension $31$.
\end{ex}

\begin{ex}\label{ex:branched-so}
For the orthogonal Lie algebras, the 
smallest example are given by $d=(1,1,2,2,1,1)$, 
i.e. (a)-type of $\liea{g}=\liea{so}_8$
and by $d=(2,2,1,2,2)$  
for an odd number of blocks 
in $\liea{so}_9$. The following branched 
diagrams give Richardson elements for the 
corresponding parabolic subalgebras. 
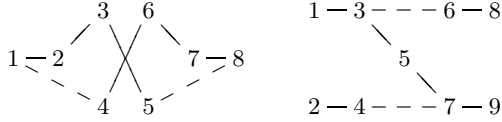
\[
{\small
\xymatrix@-6mm{
 & & 3\ar@{-}[rdd] & 6\ar@{-}[rd] \\
1\ar@{-}[r]\ar@{--}[rrd]&2\ar@{-}[ru] & & &7\ar@{-}[r]&8\\
 & & 4\ar@{-}[ruu]& 5\ar@{--}[rru]
}\quad\quad
\xymatrix@-6mm{
1\ar@{-}[r] & 3\ar@{-}[rd]\ar@{--}[rr]   & & 6\ar@{-}[r] & 8 \\
  &   & 5\ar@{-}[rd] & & \\
2\ar@{-}[r] & 4\ar@{--}[rr] & & 7\ar@{-}[r] & 9
}}
\]
\end{ex}

\begin{proof}[Proof of Lemma~\ref{lm:dual-part}]
We prove the statement for the symplectic 
Lie algebras. The corresponding statements 
for $\liea{so}_N$ are proven similarly. 

\noindent
\underline{(i) - Type (a) of $\liea{sp}_{2n}$}: 

\noindent 
Note that the bottom-top ordering of the second 
half of $L_{even}(d)$ 
ensures that the counterpart of a line 
$i$---$j$ (for $j\le n$) is again 
horizontal and that all lines 
connecting any entry of column $r$ to an 
entry to its right are horizontal. Therefore the 
line diagram $L_{even}$ has the same shape 
as the horizontal line diagram defined for 
$\liea{sl}_n$. 
In particular, the orbit of the nilpotent 
element defined by $L_{even}(d)$ 
has the same partition as the one defined 
by $L_h(d)$. Then the assertion follows with 
Lemma~\ref{lm:diagr-Jordan}.  

\noindent
\underline{(ii) - Type (b) of $\liea{sp}_{2n}$}:\\
The proof is done by induction on $r$. 
Let $d=(d_1,d_2,d_1)$ be the dimension vector. 
If $d_1\notin D_o$ (i.e. $d_1$ is not an odd entry smaller 
than $d_2$) then the line diagram $L_{even}(d_1,d_2,d_1)$ 
has the same 
chains of lines as the horizontal diagram for $\liea{sl}_{2n}$. 
For $d_1\in D_o$ the diagram $L_{even}(d_1,d_2,d_1)$ 
has $d_1-1$ two-chains 
(chains of length two) and $2$ one-chains (i.e. lines). 
So the kernel of the map $X^k$ has dimension 
$d_2$, $d_1+d_2+1$, $2d_1+d_2$ for $k=1,2,3$, 
giving the partition $1^{d_2-d_1-1}, 2^2,3^{d_1-1}$ 
and the dual of it is $d_2,d_1+1,d_1-1$ as claimed. 

Let now $d=(d_1,\dots,d_r,d_{r+1},d_r,\dots,d_1)$ with 
$d_1\le\dots\le d_{r+1}$. For $d'=(d_2,\dots,d_r,d_{r+1},d_r,
\dots,d_2)$ 
is ok. 
Let $d_1$ be even. If $d_1=d_{r+1}$ then the diagram 
$L_{odd}(d)$ is the same as $L_h(d)$, the claim follows 
immediately. 
If $d_1<d_{r+1}$, 
the diagram $L_{odd}(d)$ is obtained from $L_{odd}(d')$ 
by extending $d_1$ $(2r-2)$-chains 
to $2r$-chains. 
The kernels of the map $X^k$ satisfy 
$\dim\ker X^k=\dim\ker Y^k$ for $k\le 2r-1$, 
$\dim\ker X^{2r}=2n-d_1=\dim\ker Y^{2r}+d_1$ 
and $\dim\ker X^{2r+1}=2n=\dim\ker Y^{2r+1}+2d_1$ 
where $Y\in\liea{sp}_{2n-2d_1}$ is defined by the 
line diagram $L_{even}(d')$.  
If the partition of $Y$ is $1^{b_1},2^{b_2},
\dots,(2r-1)^{b_{2r-1}}$ 
then the partition of $X$ is 
\[
1^{b_1},\dots,(2r-2)^{b_{2r-2}}, (2r-1)^{b_{2r-1}-d_1},
(2r)^0,(2r+1)^{d_1}.
\]
Thus the dual of this partition is the dual of the partition 
of $Y$ together 
with the parts $d_1,d_1$.

If $d_1$ is even and $d_1>d_{r+1}$, the diagram $L_{odd}(d)$ 
is obtained 
from $L_{odd}(d')$ by extending $d_{r+1}$ $(2r-2)$-chains to 
$2r$-chains 
and by extending $d_1-d_{r+1}$ $(2r-3)$-chains to $(2r-1)$-chains.
Here we get 
$\dim\ker X^k =\dim\ker Y^k$ for $k\le 2r-2$, 
$\dim\ker X^{2r-1}= \dim\ker Y^{2r-1}+d_1-d_{r+1}$, 
$\dim\ker X^{2r}= 2n-d{r+1}=\dim\ker Y^{2r}+2d_1-d_{r+1}$ 
and $\dim\ker X^{2r+1}=2n=\dim\ker Y^{2r+1}+2d_1$.
So the partition of $X$ can be calculated to be 
\[
1^{b_1},\dots,(2r-3)^{b_{2r-3}},(2r-2)^{b_{2r-2}-d_1+d_{r+1}}, 
(2r-1)^{b_{2r-1}-d_{r+1}},
(2r)^{d_1-d_{r+1}},(2r+1)^{d_{r+1}}
\]
with $b_{2r-1}=d_{r+1}$. Again, the dual of the partition of 
$X$ is obtained 
from the dual of the partition of $Y$ by adding $d_1,d_1$.

Let $d_1$ be odd and $d_1>d_{r+1}$. In particuar, there are no 
odd $d_i$ that are smaller than $d_{r+1}$. The shape of 
$L_{odd}(d)$ 
is the same as the diagram for $\liea{sl}_{2n}$ (i.e. they have 
the same chain lengths). So the dual of the partition is just the 
dimension vector and we are done. 
If $d_1<d_{r+1}$, the diagram $L_{odd}(d)$ is obtained from
$L_{odd}(d')$ by extending 
$d_1-1$ $(2r-2)$-chains to $2r$-chains 
and by extending two $(2r-2)$-chains to $(2r-1)$-chains. 
The calculations 
of the dimensions of the kernels for $X$ (compared to those for 
$Y$) 
give as partition of $X$: 
\[
1^{b_1},\dots,(2r-2)^{b_{2r-2}}, (2r-1)^{b_{2r-1}-d_1-1},
(2r)^{2},(2r+1)^{d_1-1}
\]
Hence the dual of the partition of $X$ is 
obtained from 
the dual of the partition of $Y$ by adjoining $d_1+1,d_1-1$.
\end{proof}


\end{document}